\documentclass[11pt]{article}
\usepackage[rmargin=0.35in,right=0in, outer=0in,lmargin=0.35in,left=0in, inner=0in,tmargin=0.35in,bmargin=0.5in,paperwidth=5.5in,paperheight=8.5in]{geometry}
\usepackage{amsmath,amsxtra,amssymb,latexsym, amscd,amsthm}
\usepackage{mathptmx}
\usepackage{enumerate}  
\usepackage{setspace}%

\usepackage{mathtools, graphicx}
\usepackage[colorlinks=true,linkcolor=black,anchorcolor=black,citecolor=black,filecolor=black,menucolor=black,runcolor=black,urlcolor=black]{hyperref}
\usepackage[mathscr]{euscript}

\usepackage{nag, todonotes}

\newcommand{\R}{\mathbb{R}}

\newcommand{\Z}{\mathbb{Z}}

\newtheorem{theorem}{Theorem}[section]
\newtheorem{definition}{Definition}[section]

\newtheorem{lemma}[theorem]{Lemma}
\theoremstyle{remark}
\newtheorem{example}[theorem]{Example}
\newtheorem{remark}{Remark}[theorem]
\begin{document}

\title{\textbf{{\Large Characterization of Colorings Obtained by a Method of Szlam\thanks{This work was supported by NSF DMS grant no. 1950563}}}}

\author{
\small Eric Myzelev\\
\small University of Pennsylvania\\
\small myzelev@sas.upenn.edu
}

\date{}
\maketitle 
\abstract {Szlam's Lemma began life as a way of getting upper bounds on the chromatic numbers of distance graphs in normed vector spaces. Now analogs are available in a variety of hypergraph settings, but the method always involves a shrewdly chosen 2-coloring of the vertex set of a hypergraph, together with a subset of the vertex set which satisfies certain requirements with reference to the 2-coloring. From these ingredients a proper coloring of the hypergraphs is cooked up. 

In this paper, we separate the process from the conclusion of Szlam's Lemma by defining Szlam colorings of the vector spaces $\R^d$, and then a more regimented variety of these, which we call ordered Szlam colorings, which we characterize.}\\

\section{Introduction}

In \cite{2} the following strange result appeared with a proof so beautiful that Erdős very likely considered it to be From the Book. 

\begin{theorem}
Suppose that $\R^2$ is colored with red and blue so that the Euclidean distance 1 is forbidden for blue. Then the red set contains a translate of each 3-point set in $\R^2$.
\end{theorem}

Not long after, Rozalía Juhász proved \cite{6}:

\begin{theorem}
Suppose that $\R^2$ is colored with red and blue so that the Euclidean distance 1 is forbidden for blue. Then for each $F\subseteq\R^2$, $|F|=4$, the red set contains a set congruent to $F$.
\end{theorem}

This result reigned as one of the summits of Euclidean Ramsey theory for almost 40 years.

In 1999 Arthur Szlam, an undergraduate at Emory University at the time, discovered, while participating in a summer Research Experience for Undergraduates at Auburn University, a proof of Theorem 1.1 different from that in \cite{2}. After a few weeks of examining his reasoning, Arthur came upon a curious result with a very easy proof that immediately rendered Theorem 1.1 almost trivial and eventually, with the assistance of the break-through of Aubrey de Grey in \cite{3}, provided a very short proof of a strengthening of Theorem 1.2. We call his result, and each of its descendants, Szlam's Lemma.

We give here the statement and proof of an early version of Szlam's Lemma that inspired the question addressed in this paper. This version is slightly more general than the original \cite{9}, which has been greatly analogized and generalized; see \cite{1}, \cite{4}, \cite{5}, \cite{7}, and \cite{8}. We have no doubt that the question answered here can be asked and sometimes usefully answered in different settings - although in the case of the recent generalization of Szlam's Lemma in \cite{8}, even formulating the question seems to pose difficulties.

If $||\cdot||$ is a norm in $\R^d$, let $((\R^d, ||\cdot||), 1)$ stand for the unit distance graph on $(\R^d, ||\cdot ||)$ -- the graph with vertex set $\R^d$ in which $u,v\in \R^d$ are adjacent if and only if $||u-v||=1$.

\begin{lemma}[Szlam's Lemma]
Suppose that $d$ is a positive integer, $||\cdot||$ is a norm on $\R^d$, $R, B$ is a partition of $\R^d$ such that if $u,v\in B$ then $||u-v||\neq 1$, and $F\subseteq \R^d$ is such that no translate of $F$ in $\R^d$ is contained in $R$. Then $\chi((\R^d, ||\cdot||), 1)\leq |F|$.
\end{lemma}

\begin{proof}
For each $v\in R^d$, $v+F$ contains at least one point in $B$. Color $v$ with some $f\in F$ such that $v+f\in B$.

Suppose $v,w\in\R^d$ bear the same color $f\in F$ under this coloring. Then $v+f, w+f\in B$. \\
$\implies 1\neq ||(v+f)-(w-f)||=||v-w||$ \\
$\implies$ $v,w$ are not adjacent in $((\R^d, ||\cdot||), 1)$. \\
Therefore, this coloring of $((\R^d, ||\cdot||), 1)$ with $|F|$ colors is proper.
\end{proof}

When $d=2$, and $||\cdot||$ is the usual Euclidean norm on $\R^2$, Theorem 1.1 instantly follows from the Lemma and the well-known (by 1970) fact that $\chi((\R^2, ||\cdot||), 1)\geq 4$. In fact, Szlam's Lemma gives us Theorem 1.1 with Euclidean distance replaced by any distance defined by a norm $||\cdot||$ such that $\chi((\R^2, ||\cdot||), 1)\geq 4$. (This inequality holds for any of the usual norms on $||\cdot||_p$ on $\R^2$, $1\leq p\leq \infty$, with equality when $p\in \{1,\infty\}$. We suspect that it holds for all norms on $\R^2$.)

Thanks to de Grey's proof that $\chi((\R^2, ||\cdot||), 1)\geq 5$ when $||\cdot||$ is the Euclidean norm on $\R^2$ \cite{3}, from the Lemma, we have the improvement of Theorem 1.2 obtained by replacing ``set congruent to" by ``translate of" in the theorem's statement. Not only that: we are assured that the improvement holds when the Euclidean norm is replaced by any norm $||\cdot||$ such that $\chi((\R^2, ||\cdot||), 1)> 4$.

\section{The Question}

Suppose that $R, B$ is a partition of $\R^d$  and $F\subseteq \R^d$ is a non-empty set no translate of which is contained in $R$.

\begin{definition}
A \underline{Szlam coloring} of $\R^d$ associated with $R,B$ and $F$ is a function $\varphi:\R^d\to F$ such that for each $v\in \R^d$, $v+\varphi(v)\in B$.
\end{definition}
\begin{definition}
If $|F|<\infty$ and $f_1,...,f_k$ is an ordering of $F$, the \underline{ordered Szlam coloring} of $\R^d$ associated with $R,B$ and $F$ and the given ordering of $F$ is the function $\varphi:\R^d\to \{1,...,k\}$ such that for each $v\in \R^d$, $\varphi(v)$ is the smallest index $i$ such that $v+f_i\in B$.
\end{definition}

By definition, any coloring produced in an application of Lemma 1.3 (Szlam's Lemma) will in fact be a Szlam coloring. However, it is possible that for some $v\in\R^d$ there exist multiple elements $f\in F$ such that $v+f\in B$, in which case we choose between these arbitrarily. It was our ambition to find a useful characterization of Szlam colorings, but due to the potentially uncountable number of arbitrary decisions that are made, now we are inclined to the opinion that the very existence of such a characterization may be a question suitable for philosophical debate. 

We can, however, tame the set of possible colorings if we fix an ordering on the vertices of $F$, and always select the first vertex in $F$ when given a choice. An application of Lemma 1.3 with a fixed ordering always yields an ordered Szlam coloring. In this paper we provide a characterization of ordered Szlam colorings.

\section{The Answer}

Suppose that $\varphi:\R^d\to C$, $|C|=k>1$, is a coloring of $\R^d$. Let $c_1,...,c_k$ be an ordering of $C$, and let $A_i=\varphi^{-1}(\{c_i\})$, $i=1,...,k$. We will say that $\varphi$ is dominant with respect to the ordering $c_1,...,c_k$ of $C$ if and only if there exists distinct vectors $t_2,...,t_k\subset R^d$ such that
\begin{itemize}
\item $A_i+t_i\subseteq A_1$ for all $i>1$, and
\item for all $1<i<j$, $(A_j+t_i)\cap A_1=\emptyset$.
\end{itemize}
That is, for $i\in\{2,...,k\}$, translation by $t_i$ sends $A_i$ into $A_1$, and for each $j$ such that $i<j$, translation by $t_i$ sends $A_j$ into $\R^d\setminus A_1$.

\begin{theorem}
A coloring $\varphi:\R^d\to C$, $|C|=k$, is an ordered Szlam coloring if and only if there is an ordering $c_1,...,c_k$ of $C$ with respect to which $\varphi$ is dominant. Further, for every such ordering of $C$, $\varphi$ is an ordered Szlam coloring of $\R^d$ with reference to $R, B, F,$ and an ordering $f_1,...,f_k$ of $F$ which corresponds to the ordering of $C$: for $j\in\{1,...,k\}$, and $v\in \R^d$, $\varphi(v)=c_j\Leftrightarrow v+f_j\in B$ and for all $i<j$, $v+f_i\in R$.
\end{theorem}

\begin{proof}
Suppose that $R,B$ partition $\R^d$, $F\subseteq \R^d$, no translate of $F$ is contained in $R$, and $f_1,...,f_k$ is an ordering of $F$. Let $\varphi:\R^d\to\{1,...,k\}$ be the ordered Szlam coloring associated with $R,B,F$, and the given ordering of $F$. For $j\in\{1,...,k\}$ let $A_j=\varphi^{-1}(\{j\})$.

If $a\in B-f_1$ then $a+f_1\in B$, so $\varphi(a)=1$. Therefore $B-f_1\subseteq A_1$. On the other hand, $a\in A_1\implies \varphi(a)=1\implies a+f_1\in B\implies a\in B-f_1$. Thus $A_1=B-f_1$.

Now suppose that $1<j\leq k$. For $a\in A_j$, $\varphi(a)=j\implies a+f_j\in B$ and for all $i\in\{1,...,j-1\}$, $a+f_i\in R=\R^d\setminus B$. Thus $A_j\subseteq B-f_j=A_1+(f_1-f_j)$, while for $i<j$, $(A_j+f_i)\cap B=\emptyset\implies (A_j+f_i-f_1)\cap (B-f_1)=(A_j+f_i-f_1)\cap A_1=\emptyset$.

Thus $\varphi$ is dominant with respect to the ordering $1,...,k$ of the color set, with translaters $t_j=f_j-f_1$, $j\in\{2,...k\}$.

Now suppose that $\varphi:\R^d\to C$ is a dominant coloring of $\R^d$ with respect to an ordering $c_1,...,c_k$ of $C$, with distinct translaters $t_2,...,t_k$, such that, with $A_i=\varphi^{-1}(\{c_i\})$, $i=1,...,k$, $A_i+t_i\subseteq A_1$, for all $i\geq 2$ and, for all $2\leq i<j\leq k$, $(A_j+t_i)\cap A_1=\emptyset$. Set $B=A_1, R=\R^d\setminus B$, and $F=\{f_1,...,f_k\}$, with ordering $f_1,...,f_k$, with $f_1=0$, $f_i=t_i, i=2,...,k$. 

Observe that for any $v\in\R^d$, if $\varphi(v)=c_1$ then $v+f_1=v\in B$, and if $\varphi(v)=c_j, j>1$ then $v+f_j\in A_j+t_j\subseteq A_1=B$, so $R, B, F$, and the ordering $f_1,...,f_k$ of $F$ are allowable as a foundation for an ordered Szlam coloring of $\R^d$. Let $\psi$ be that coloring. We will show that $\varphi=\psi$, in the sense that $\varphi(v)=c_i\Leftrightarrow \psi(v)=i$.

If $\varphi(v)=c_1$ then $v\in A_1\implies v\in 0+A_1=f_1+A_1\implies \psi(v)=1$. Otherwise, if $\varphi(v)=c_j, j>1$, then $v\in A_j\implies v+f_j\in A_1$, and $v+f_i\notin A_1$ for $i<j\implies \varphi(v)=j$.
\end{proof}

\begin{example}
Suppose that $\R$ is colored with green, yellow, and white as follows: $\varphi(r)=$ green if $r\in [6k, 6k+3)$, for some $k\in\Z$, $\varphi(r)=$ yellow if $r\in [6k+3, 6k+5)$ for some $k\in\Z$, $\varphi(r)=$ white if $r\in [6k+5, 6k+6)$ for some $k\in\Z$.

Letting $G$ denote the set of green points, $Y$ the set of yellow points, and $W$ the set of white points, we ask if this coloring is an ordered Szlam coloring of $\R$, and, if so, what are the constituents $R, B,$ and $F=\{f_1,f_2,f_3\}$ with reference to which $\varphi$ is defined?

Clearly
$$G=\bigcup_{k\in\Z}[6k, 6k+3)$$
$$Y = \bigcup_{k\in\Z} [6k+3, 6k+5)$$
$$W = \bigcup_{k\in\Z} [6k+5, 6k+6)$$

Clearly neither $Y$ nor $W$ can contain a translate of $G$. Therefore, we take $G=A_1$ and green $=c_1$, and search for $t_2,t_3$ and a choice of yellow, white $\in\{c_2,c_3\}$.

We take yellow $=c_2$, white $=c_3$, and $t_2=-2$, $t_3=-3$ to see that this coloring is, indeed, an ordered Szlam coloring. Back in Szlam-land, $B=\bigcup_{k\in \Z}[6k, 6k+3), R=\R\setminus B$, and $F=\{0,-2,-3\}$ with ordering $f_1=0,f_2=-2,f_3=-3$.
\end{example}

\begin{remark}
If $F=\{0,-2,-3\}$ is replaced by a translate of $F$, say $\tilde{F}=\{r, r-2,r-3\}$, which will satisfy the requirement that no translate of $\tilde{F}$ is contained in $R$, and if the ordering of $\tilde{F}$ is $\tilde{f_1}=r$, $\tilde{f_2}=r-2,\tilde{f_3}=r-3$, the ordered Szlam coloring defined will not be the $\psi=\varphi$ we got before: it will be that coloring "translated by r", meaning 
$$\tilde{G}=\bigcup_{k\in\Z}[6k-r, 6k-r+3)=G-r$$
$$\tilde{Y}=Y-r, \tilde{W}=W-r$$
\end{remark}

\begin{remark}
In the Hadwiger-Isbell proof that $\chi(\R^2,1)\leq 7$, the plane is tiled with regular hexagons of diameter $d$, $\frac{2}{\sqrt{7}}<d\leq 1$. These hexagons are grouped into "Hadwiger tiles," 7 hexagons consisting of one central hexagons together with the 6 hexagons surrounding it. Color one Hadwiger tile with 7 colors, one to a hexagon, and copy this coloring periodically on all the other Hadwiger tiles. If $\frac{2}{\sqrt{7}}<d<1$ - and the aim is to forbid Euclidean distance 1 for each of the 7 colors, then no care need be taken about coloring hexagon boundaries.

However, as shown in \cite{7}, if one does take care about coloring the boundaries, then each of these colorings can be transformed into a Szlam coloring. The set $F$ is the set of hexagon centers in one particular tile, and the set $B$ is the union of the closures of the central hexagons in the Hadwiger tiles. The moral of this story is that in using Szlam's Lemma to put upper bounds on the chromatic numbers of distance graphs on $(\R^d, ||\cdot ||)$, you may have to work harder than absolutely necessary. 

Still there is left standing a question of interest independent of which of 5,6,7 equals $\chi(\R^2,1)$: whatever $k=\chi(\R^2,1)$ may be, is there a Szlam coloring $(R,B,F)$ of $\R^2$ such that $B$ contains no two points a distance 1 apart and $|F|=k$? Of course, if $k=7$, then, by \cite{7}, the answer is yes.
\end{remark}

\section*{Acknowledgement}

We would like to acknowledge that the first draft of this work was written with much advice from Dr. Peter Johnson, the author's mentor during the 2023 NSF-funded summer REU program at Auburn University. The author extends his deepest gratitude for Dr. Johnson's generous encouragement and constructive criticism throughout the program.

\end{document}